\documentclass[12pt]{article}%
\usepackage{amsfonts}
\usepackage{makeidx}
\usepackage{amsmath,exscale,amssymb}
\usepackage{float}
\usepackage{graphicx}
\usepackage{amscd}
\usepackage{amsthm, amsxtra}
\usepackage{times}
\usepackage{verbatim}
\usepackage{amsmath}
\usepackage{amssymb}%
\theoremstyle{plain}
\newtheorem{X}{X}
\newtheorem*{theorem}{Theorem}

\newtheorem{corollary}[X]{Corollary}

\theoremstyle{definition}

\newcounter{junk}

\newenvironment{mylist}
{\begin{list}
{--}
{\setlength{\leftmargin}{.5in}\setlength{\rightmargin}{.5in}}}
{\end{list}}

\def\1{{1\mkern-7mu1}}  

\newcommand\Gal{\operatorname{Gal}}

\newcommand\Hom{\operatorname{Hom}}

\newcommand\Ker{\operatorname{Ker}}

\begin{document}

\title{The de Rham-Witt and ${\mathbb{Z}}_{p}$-cohomologies of an algebraic variety}
\author{James S. Milne
\and Niranjan Ramachandran\thanks{Partially supported by GRB (University of
Maryland) and IHES.}}
\date{\today}
\maketitle

\hfill\textit{To Mike Artin on the occasion of his 70th birthday.}

\begin{abstract}
We prove that, for a smooth complete variety $X$ over a perfect field,
\[
H^{i}(X,\mathbb{Z}{}_{p}(r))\cong\Hom_{\mathsf{D}{}_{c}^{b}(R)}(\1,R\Gamma
(W\Omega_{X}^{\bullet})(r)[i])
\]
where $H^{i}(X,\mathbb{Z}{}_{p}(r))=\varprojlim_{n}H^{i-r}(X_{\text{et}}%
,\nu_{n}(r))$ (Milne 1986, p309), $W\Omega_{X}^{\bullet}$ is the de Rham-Witt
complex on $X$ (Illusie 1979b), and $\mathsf{D}_{c}^{b}(R)$ is the
triangulated category of coherent complexes over the Raynaud ring (Illusie and
Raynaud 1983, I 3.10.1, p120).

\end{abstract}

\subsection{Introduction}

According to the standard philosophy (cf. Deligne 1994, 3.1), a cohomology
theory $X\mapsto H^{i}(X,r)$ on the algebraic varieties over a fixed field $k$
should arise from a functor $R\Gamma$ taking values in a triangulated category
$\mathsf{D}$ equipped with a $t$-structure and a Tate twist $D\mapsto D(r)$ (a
self-equivalence). The heart $\mathsf{D}{}^{\heartsuit}$ of $\mathsf{D}$
should be stable under the Tate twist and have a tensor structure; in
particular, there should be an essentially unique identity object $\1$ in
$\mathsf{D}^{\heartsuit}$ such that $\1\otimes D\cong D\cong D\otimes\1$ for
all objects in $\mathsf{D}^{\heartsuit}$. The cohomology theory should satisfy%
\begin{equation}
H^{i}(X,r)\cong\Hom_{\mathsf{D{}}}(\1,R\Gamma(X)(r)[i]). \label{e2}%
\end{equation}
For example, motivic cohomology $H_{\mathrm{mot}}^{i}(X,\mathbb{Q}{}(r))$
should arise in this way from a functor to a category $\mathsf{D{}}$ whose
heart is the category of mixed motives $k$. Absolute $\ell$-adic \'{e}tale
cohomology $H_{\mathrm{et}}^{i}(X,\mathbb{Z}_{\ell}(r))$, $\ell\neq
\mathrm{char}(k)$, arises in this way from a functor to a category
$\mathsf{D}$ whose heart is the category of continuous representations of
$\Gal(\bar{k}/k)$ on finitely generated $\mathbb{Z}{}_{\ell}$-modules (Ekedahl
1990). When $k$ is algebraically closed, $H_{\mathrm{et}}^{i}(X,\mathbb{Z}%
{}_{\ell}(r))$ becomes the familiar group $\varprojlim H_{\mathrm{et}}%
^{i}(X,\mu_{\ell^{n}}^{\otimes r})$ and lies in $\mathsf{D}^{\heartsuit}$;
moreover, in this case, (\ref{e2}) simplifies to
\begin{equation}
H^{i}(X,r)\cong H^{i}(R\Gamma(X)(r)). \label{e1}%
\end{equation}

Now let $k$ be a perfect field of characteristic $p\neq0$, and let $W$ be the
ring of Witt vectors over $k$. For a smooth complete variety $X$ over $k$, let
$W\Omega_{X}^{\bullet}$ denote the de Rham-Witt complex of
Bloch-Deligne-Illusie (see Illusie 1979b). Regard $\Gamma=\Gamma(X,-)$ as a
functor from sheaves of $W$-modules on $X$ to $W$-modules. Then%
\[
H_{\mathrm{crys}}^{i}(X/W)\cong H^{i}(R\Gamma(W\Omega_{X}^{\bullet}))
\]
(Illusie 1979a, 3.4.3), where $H_{\mathrm{crys}}^{i}(X/W)$ is the crystalline
cohomology of $X$ (Berthelot 1974). In other words, $X\mapsto H_{\mathrm{crys}%
}^{i}(X/W)$ arises as in (\ref{e1}) from the functor $X\mapsto R\Gamma
(W\Omega_{X}^{\bullet})$ with values in $\mathsf{D{}}^{+}(W)$.

Let $R$ be the Raynaud ring, let $\mathsf{D}(X,R)$ be the derived category of
the category of sheaves of graded $R$-modules on $X$, and let $\mathsf{D}(R)$
be the derived category of the category of graded $R$-modules (Illusie 1983,
2.1). Then $\Gamma$ derives to a functor
\[
R\Gamma\colon\mathsf{D}(X,R)\rightarrow\mathsf{D}(R).
\]
When we regard $W\Omega_{X}^{\bullet}$ as a sheaf of graded $R$-modules on
$X$, $R\Gamma(W\Omega_{X}^{\bullet})$ lies in the full subcategory
$\mathsf{D}_{c}^{b}(R)$ of $\mathsf{D}(R)$ consisting of coherent complexes
(Illusie and Raynaud 1983, II 2.2), which Ekedahl has shown to be a
triangulated subcategory with $t$-structure (Illusie 1983, 2.4.8). In this
note, we define a Tate twist $(r)$ on $\mathsf{D{}}_{c}^{b}(R)$ and prove that%
\[
H^{i}(X,\mathbb{Z}{}_{p}(r))\cong\Hom_{\mathsf{D{}}_{c}^{b}(R)}(\1,R\Gamma
(W\Omega_{X}^{\bullet})(r)[i]).
\]
Here $H^{i}(X,\mathbb{Z}{}_{p}(r))=_{\mathrm{df}}\varprojlim_{n}%
H_{\mathrm{et}}^{i-r}(X,\nu_{n}(r))$ with $\nu_{n}(r)$ the additive subsheaf
of $W_{n}\Omega_{X}^{r}$ locally generated for the \'{e}tale topology by the
logarithmic differentials (Milne 1986, \S 1), and $\1$ is the identity object
for the tensor structure on graded $R$-modules defined by Ekedahl (Illusie
1983, 2.6.1). In other words, $X\mapsto H^{i}(X,\mathbb{Z}{}_{p}(r))$ arises
as in (\ref{e2}) from the functor $X\mapsto R\Gamma(W\Omega_{X}^{\bullet})$
with values in $\mathsf{D}_{c}^{b}(R)$.

This result is used in the construction of the triangulated category of
integral motives in Milne and Ramachandran 2005.

It is a pleasure for us to be able to contribute to this volume: the
$\mathbb{Z}{}_{p}$-cohomology was introduced (in primitive form) by the first
author in an article whose main purpose was to prove a conjecture of Artin,
and, for the second author, Artin's famous 18.701-2 course was his first
introduction to real mathematics.

\subsection*{The Tate twist}

According to the standard philosophy, the Tate twist on motives should be
$N\mapsto N(r)=N\otimes\mathbb{T}{}^{\otimes r}$ with $\mathbb{T}{}$ dual to
$\mathbb{L}{}$ and $\mathbb{L}{}$ defined by $Rh({\mathbb{P}}^{1}%
)=\1\oplus{\mathbb{L}}[-2]$.

The Raynaud ring is the graded $W$-algebra $R=R^{0}\oplus R^{1}$ generated by
$F$ and $V$ in degree $0$ and $d$ in degree $1$, subject to the relations
$FV=p=VF$, $Fa=\sigma a\cdot F$, $aV=V\cdot\sigma a$, $ad=da$ ($a\in W$),
$d^{2}=0$, and $FdV=d$; in particular, $R^{0}$ is the Dieudonn\'{e} ring
$W_{\sigma}[F,V]$ (Illusie 1983, 2.1). A graded $R$-module is nothing more
than a complex
\[
M^{\bullet}=(\cdots\rightarrow M^{i}\overset{d}{\longrightarrow}%
M^{i+1}\rightarrow\cdots)
\]
of $W$-modules whose components $M^{i}$ are modules over $R^{0}$ and whose
differentials $d$ satisfy $FdV=d$. We define $T$ to be the functor of graded
$R$-modules such that $(TM)^{i}=M^{i+1}$ and $T(d)=-d$. It is exact and
defines a self-equivalence $T\colon\mathsf{D}_{c}^{b}(R)\rightarrow
\mathsf{D}_{c}^{b}(R)$.

The identity object for Ekedahl's tensor structure on the graded $R$-modules
is the graded $R$-module
\[
\1=(W,F=\sigma,V=p\sigma^{-1})
\]
concentrated in degree zero (Illusie 1983, 2.6.1.3). It is equal to the module
$E_{0/1}=_{\mathrm{df}}R^{0}/(F-1)$ of Ekedahl 1985, p. 66.

There is a canonical homomorphism%
\[
{\1}{\oplus}T^{{-1}}({\1})[-1]\rightarrow R\Gamma(W\Omega_{\mathbb{P}{}^{1}%
}^{\bullet})
\]
(in $\mathsf{D}_{c}^{b}(R)$), which is an isomorphism because it is on
$W_{1}\Omega_{\mathbb{P}{}^{1}}^{\bullet}=\Omega_{\mathbb{P}{}^{1}}^{\bullet}$
and we can apply Ekedahl's \textquotedblleft Nakayama lemma\textquotedblright%
\ (Illusie 1983, 2.3.7). See Gros 1985, I 4.1.11, p21, for a more general
statement. This suggests our definition of the Tate twist $r$ (for $r\geq0$),
namely, we set
\[
M{}(r)=T^{r}(M)[-r]
\]
for $M$ in $\mathsf{D}_{c}^{b}(R)$.

Ekedahl has defined a nonstandard $t$-structure on $\mathsf{D}_{c}^{b}(R)$ the
objects of whose heart $\Delta$ are called diagonal complexes (Illusie 1983,
6.4). It will be important for our future work to note that $\mathbb{T=}%
{}T({\1})[-1]$ is a diagonal complex: the sum of its module degree $(-1)$ and
complex degree $(+1)$ is zero. The Tate twist is an exact functor which
defines a self-equivalence of $\mathsf{D}_{c}^{b}(R)$ preserving $\Delta$.

\subsection{Theorem and corollaries}

Regard $W\Omega_{X}^{\bullet}$ as a sheaf of graded $R$-modules on $X$, and
write $R\Gamma$ for the functor $\mathsf{D}(X,R)\rightarrow\mathsf{D}(R)$
defined by $\Gamma(X,-)$. As we noted above, $R\Gamma(W\Omega_{X}^{\bullet})$
lies in $\mathsf{D}_{c}^{b}(R)$.

\begin{theorem}
For any smooth complete variety $X$ over a perfect field $k$ of characteristic
$p\neq0$, there is a canonical isomorphism%
\[
H^{i}(X{},\mathbb{Z}{}_{p}(r))\cong\Hom_{\mathsf{D}_{c}^{b}(R)}(\1,R\Gamma
(W\Omega_{X}^{\bullet}){}(r)[i]).
\]

\end{theorem}

\begin{proof}
For a graded $R$-module $M^{\bullet}$,
\[
\Hom(\1,M^{\bullet})=\Ker(1-F\colon M^{0}\rightarrow M^{0}).
\]
To obtain a similar expression in $\mathsf{D}^{b}(R)$ we argue as in Ekedahl
1985, p90. Let $\hat{R}$ denote the completion $\varprojlim R/(V^{n}%
R+dV^{n}R)$ of $R$ (ibid. p60). Then right multiplication by $1-F$ is
injective, and $\1\cong\hat{R}^{0}/\hat{R}^{0}(1-F)$. As $F$ is topologically
nilpotent on $\hat{R}^{1}$, this shows that the sequence%
\begin{equation}
\begin{CD} 0@>>>\hat{R}@>{\cdot(1-F)}>>\hat{R}@>>>\1@>>>0,\end{CD} \label{e3}%
\end{equation}
is exact. Thus, for a complex of graded $R$-modules $M$ in $\mathsf{D}^{b}%
(R)$,%
\[
\Hom_{\mathsf{D}(R)}(\1,M)\overset{\text{Grivel 1987, 10.9}}{\cong}%
H^{0}(R\Hom(\1,M))\overset{(\ref{e3})}{\cong}H^{0}(R\Hom(\hat{R}%
\xrightarrow{\cdot(1-F)}\hat{R},M)).
\]
If $M$ is complete in the sense of Illusie 1983, 2.4, then $R\Hom(\hat
{R},M)\cong R\Hom(R,M)$ (Ekedahl 1985, I 5.9.3ii, p78), and so%
\begin{align}
\Hom_{\mathsf{D}(R)}(\1,M)  &  \cong H^{0}%
(\Hom(R\xrightarrow{\cdot(1-F)}R,M))\nonumber\\
&  \cong H^{0}(\Hom(R,M)\xrightarrow{1-F}\Hom(R,M)). \label{e9}%
\end{align}

Following Illusie 1983, 2.1, we shall view a complex of graded $R$-modules as
a bicomplex $M^{\bullet\bullet}$ in which the first index corresponds to the
$R$-grading: thus the $j^{\mathrm{th}}$ row $M^{\bullet j}$ of the bicomplex
is the $R$-module $(\cdots\rightarrow M^{i,j}\rightarrow M^{i+1,j}%
\rightarrow\cdots)$, and the $i^{\mathrm{th}}$ column $M^{i\bullet}$ is a
complex of (ungraded) $R^{0}$-modules. The $j^{\text{th}}$-cohomology
$H^{j}(M^{\bullet\bullet})$ of $M^{\bullet\bullet}$ is the graded $R$-module
\[
(\cdots\rightarrow H^{j}(M^{i\bullet})\rightarrow H^{j}(M^{i+1\bullet
})\rightarrow\cdots).
\]
Now, $\Hom(R,M^{\bullet\bullet})=M^{0\bullet}$, and so%
\begin{equation}
H^{0}(\Hom(R,M^{\bullet\bullet}(r)[i]))=H^{i-r}(M^{r\bullet}). \label{e6a}%
\end{equation}

The complex of graded $R$-modules $R\Gamma(W\Omega_{X}^{\bullet})$ is complete
(Illusie 1983, 2.4, Example (b), p33), and so (\ref{e9}) gives an isomorphism%
\begin{equation}%
\begin{split}
\Hom_{\mathsf{D}(R)}  &  (\1,R\Gamma(W\Omega_{X}^{\bullet})(r)[i])\cong\\
&  H^{0}(\Hom(R,R\Gamma(W\Omega_{X}^{\bullet}%
)(r)[i])\xrightarrow{1-F}\Hom(R,R\Gamma(W\Omega_{X}^{\bullet})(r)[i])).
\end{split}
\label{e10}%
\end{equation}
The $j^{\mathrm{th}}$-cohomology of $R\Gamma(W\Omega_{X}^{\bullet})$ is
obviously%
\[
H^{j}(R\Gamma(W\Omega_{X}^{\bullet}))=(\cdots\rightarrow H^{j}(X,W\Omega
_{X}^{i})\rightarrow H^{j}(X,W\Omega_{X}^{i+1})\rightarrow\cdots)
\]
(Illusie 1983, 2.2.1), and so (\ref{e6a}) allows us to rewrite (\ref{e10}) as
\[
\Hom_{\mathsf{D}(R)}(\1,R\Gamma(W\Omega_{X}^{\bullet})(r)[i])\cong
H^{i-r}(R\Gamma(W\Omega_{X}^{r})\xrightarrow{1-F}R\Gamma(W\Omega_{X}^{r})).
\]
This gives an exact sequence%
\begin{equation}
\cdots\rightarrow\Hom(\1,R\Gamma(W\Omega_{X}^{\bullet}){}(r)[i])\rightarrow
H^{i-r}(X{},W\Omega_{X{}}^{r})\xrightarrow{1-F}H^{i-r}(X{},W\Omega_{X{}}%
^{r})\rightarrow\cdots\label{e4}%
\end{equation}

On the other hand, there is an exact sequence (Illusie 1979b, I 5.7.2)%
\[
0\rightarrow\nu_{\bullet}(r)\rightarrow W_{\bullet}\Omega_{X}^{r}%
\xrightarrow{1-F}W_{\bullet}\Omega_{X}^{r}\rightarrow0
\]
of prosheaves on $X_{\mathrm{et}}$, which gives rise to an exact sequence%
\begin{equation}
\cdots\rightarrow H^{i}(X,\mathbb{Z}{}_{p}(r))\rightarrow H^{i-r}%
(X{},W_{\bullet}\Omega_{X{}}^{r})\overset{1-F}{\longrightarrow}H^{i-r}%
(X{},W_{\bullet}\Omega_{X{}}^{r})\rightarrow\cdots\label{e5}%
\end{equation}
(Milne 1986, 1.10). Here $\nu_{\bullet}(r)$ denotes the projective system
$(\nu_{n}(r))_{n\geq0}$, and $H^{i}(X,W_{\bullet}\Omega_{X}^{r})=\varprojlim
_{n}H^{i}(X,W_{n}\Omega_{X}^{r})$ (\'{e}tale or Zariski cohomology --- they
are the same).

Since $H^{r}(X,W\Omega_{X}^{r})\cong H^{r}(X,W_{\bullet}\Omega_{X}^{r})$
(Illusie 1979a, 3.4.2, p101), the sequences (\ref{e4}) and (\ref{e5}) will
imply the theorem once we check that there is a suitable map from one sequence
to the other, but the right hand square in%
\[
\begin{CD}
W\Omega_{X}^{r} @>{1-F}>> W\Omega_{X}^{r}\\
@VVV@VVV\\
W_{\bullet}\Omega_{X}^{r} @>{1-F}>> W_{\bullet}\Omega_{X}^{r}%
\end{CD}\qquad\xrightarrow{R\Gamma}\qquad\begin{CD}
R\Gamma W\Omega_{X}^{r} @>{1-F}>> R\Gamma W\Omega_{X}^{r}\\
@VVV@VVV\\
R\Gamma W_{\bullet}\Omega_{X}^{r} @>{1-F}>> R\Gamma W_{\bullet}\Omega_{X}^{r}%
\end{CD}
\]
gives rise to such a map.
\end{proof}

As in Milne 1986, p309, we let $H^{i}(X,(\mathbb{Z}{}/p^{n}\mathbb{Z}%
{})(r))=H_{\mathrm{et}}^{i-r}(X,\nu_{n}(r)).$

\begin{corollary}
There is a canonical isomorphism%
\[
H^{i}(X,(\mathbb{Z}/p^{n}\mathbb{Z)}(r))\cong\Hom_{D_{c}^{b}(R)}(\1,R\Gamma
W_{n}\Omega_{X}^{\bullet}(r)[i]).
\]

\end{corollary}

\begin{proof}
The canonical map $\nu_{\bullet}(r)/p^{n}\nu_{\bullet}(r)\rightarrow\nu
_{n}(r)$ is an isomorphism (Illusie 1979b, I 5.7.5, p. 598), and the canonical
map $W\Omega_{X}^{\bullet}/p^{n}W\Omega_{X}^{\bullet}\rightarrow W_{n}%
\Omega_{X}^{\bullet}$ is a quasi-isomorphism (ibid. I 3.17.3, p577). The
corollary now follows from the theorem by an obvious five-lemma argument.
\end{proof}

Lichtenbaum (1984) conjectures the existence of a complex $\mathbb{Z}{}(r)$ on
$X_{\mathrm{et}}$ satisfying certain axioms and sets $H_{\mathrm{mot}}%
^{i}(X,r)=$ $H_{\mathrm{et}}^{i}(X,\mathbb{Z}{}(r))$. Milne (1988, p68) adds
the \textquotedblleft Kummer $p$-sequence\textquotedblright\ axiom that there
be an exact triangle
\[
\mathbb{Z}{}(r)\overset{p^{n}}{\longrightarrow}\mathbb{Z}{}(r)\rightarrow
\nu_{n}(r)[-r]\rightarrow\mathbb{Z}{}(r)[1].
\]
Geisser and Levine (2000, Theorem 8.5) show that the higher cycle complex of
Bloch (on $X_{\mathrm{et}}$) satisfies this last axiom, and so we have the
following result.

\begin{corollary}
Let $\mathbb{Z}{}(r)$ be the higher cycle complex of Bloch on $X_{\mathrm{et}%
}$. Then there is a canonical isomorphism%
\[
H_{\mathrm{et}}^{i}(X,\mathbb{Z}{}(r)\overset{p^{n}}{\longrightarrow
}\mathbb{Z}{}(r))\cong\Hom_{D_{c}^{b}(R)}(\1,R\Gamma W_{n}\Omega_{X}^{\bullet
}(r)[i]).
\]

\end{corollary}

\paragraph{Acknowledgement.}

We thank P. Deligne for pointing out a misstatement in the introduction to the
original version.

\subsection{References}

\small\setlength{\parindent}{0in} \setlength{\parskip}{3pt}

Berthelot, Pierre. Cohomologie cristalline des sch\'{e}mas de
caract\'{e}ristique $p>0$. Lecture Notes in Mathematics, Vol. 407.
Springer-Verlag, Berlin-New York, 1974.

Deligne, Pierre: A quoi servent les motifs? Motives (Seattle, WA, 1991),
143--161, Proc. Sympos. Pure Math., 55, Part 1, Amer. Math. Soc., Providence,
RI, 1994.

Ekedahl, Torsten: On the multiplicative properties of the de Rham-Witt
complex. II. Ark. Mat. 23, no. 1, 53--102 (1985).

Ekedahl, Torsten: Diagonal complexes and $F$-gauge structures. Travaux en
Cours. Hermann, Paris (1986).

Ekedahl, Torsten. On the adic formalism. The Grothendieck Festschrift, Vol.
II, 197--218, Progr. Math., 87, Birkh\"{a}user Boston, Boston, MA, 1990.

Geisser, Thomas; Levine, Marc: The $K$-theory of fields in characteristic $p$.
Invent. Math. 139, no. 3, 459--493 (2000).

Grivel, Pierre-Paul: Cat\'egories d\'eriv\'es et foncteurs d\'eriv\'es, in
Borel, A.; Grivel, P.-P.; Kaup, B.; Haefliger, A.; Malgrange, B.; Ehlers, F.
Algebraic $D$-modules. Perspectives in Mathematics, 2. Academic Press, Inc.,
Boston, MA, 1987.

Gros, Michel: Classes de Chern et classes de cycles en cohomologie de
Hodge-Witt logarithmique. Bull. Soc. Math. France M\'{e}m, 21, 1-87 (1985).

Illusie, Luc: Complexe de de Rham-Witt. Journ\'{e}es de G\'{e}om\'{e}trie
Alg\'{e}brique de Rennes (Rennes, 1978), Vol. I, pp. 83--112, Ast\'{e}risque,
63, Soc. Math. France, Paris, 1979a.

Illusie, Luc: Complex de de Rham-Witt et cohomologie crystalline. Ann. Scient.
\'{E}c. Norm. Sup. 12, 501--661 (1979b).

Illusie, Luc: Finiteness, duality, and K\"{u}nneth theorems in the cohomology
of the de Rham-Witt complex. Algebraic geometry (Tokyo/Kyoto, 1982), 20--72,
Lecture Notes in Math., 1016, Springer, Berlin (1983).

Illusie, Luc; Raynaud, Michel: Les suites spectrales associ\`{e}es au complexe
de de Rham-Witt. Inst. Hautes. \'{E}tudes Sci. Publ. Math. No. 57, 73--212 (1983).

Lichtenbaum, S: Values of zeta-functions at nonnegative integers. Number
theory, Noordwijkerhout 1983 (Noordwijkerhout, 1983), 127--138, Lecture Notes
in Math., 1068, Springer, Berlin, 1984.

Milne, James S: Values of zeta functions of varieties over finite fields.
Amer. J. Math. 108, no. 2, 297--360 (1986).

Milne, James S: Motivic cohomology and values of zeta functions. Compositio
Math. 68, 59-102 (1988).

Milne, James S; Ramachandran, Niranjan: The $t$-category of integral motives
and values of zeta functions. In preparation, 2005.

\bigskip

\noindent{\footnotesize James S. Milne, 2679 Bedford Rd., Ann Arbor, MI 48104,
USA, math@jmilne.org, www.jmilne.org/math/. \smallskip}

\noindent{\footnotesize Niranjan Ramachandran, Dept.\ of Mathematics,
University of Maryland, College Park, MD 20742, USA, atma@math.umd.edu,
www.math.umd.edu/$\sim$atma. }

\end{document}